\documentclass[10pt,twoside,reqno]{amsart}
\usepackage{amssymb}
\usepackage{multirow}
\calclayout

\numberwithin{equation}{section}
\makeatletter

\renewcommand{\@secnumfont}{\bfseries}

\renewcommand{\section}{\@startsection{section}{1}%
  {0mm}{.7\linespacing\@plus\linespacing}{.5\linespacing}
  {\normalfont\bfseries\centering}}

\newcommand{\bibsection}{\@startsection{section}{1}%
  {0mm}{.7\linespacing\@plus\linespacing}{.5\linespacing}
  {\normalfont\scshape\centering}}

\renewcommand{\@biblabel}[1]{#1.}

\newtheorem{thm}{\bf Theorem}[section]
\newtheorem{lem}[thm]{\bf Lemma}
\newtheorem{cor}[thm]{\bf Corollary}

\begin{document}

\vspace{1.3cm}

\title{On central complete Bell polynomials}

\author{Taekyun Kim}
\address{Department of Mathematics, Kwangwoon University, Seoul 139-701, Republic
of Korea}
\email{tkkim@kw.ac.kr}

\author{Dae San Kim}
\address{Department of Mathematics, Sogang University, Seoul 121-742, Republic
of Korea}
\email{dskim@sogang.ac.kr}

\author{Gwan-Woo Jang}
\address{Department of Mathematics, Kwnagwoon University, Seoul 139-701, Republic of Korea}
\email{gwjang@kw.ac.kr}

\maketitle

\begin{abstract}
In this paper, we consider central complete and incomplete Bell polynomials which are generalizations of the recently introduced central Bell polynomials and 'central' analogues for the complete and incomplete Bell polynomials. We investigate some properties and identities for these polynomials. Especially, we give explicit formulas for the central complete and incomplete Bell polynomials related to central factorial numbers of the second kind.
\end{abstract}
\bigskip
\medskip

\markboth{\centerline{\scriptsize On central complete Bell polynomials }}
{\centerline{\scriptsize T. Kim, D. S. Kim and G.-W. Jang}}

\section{\bf Introduction}
The Stirling numbers of the second kind are given by 
\begin{equation} \begin{split} \label{01}
\frac{1}{k!}(e^t-1)^{k}=\sum_{n=k}^\infty S_{2}(n,k)\frac{t^n}{n!},\,\,\,\,
(\textnormal{see}\,\,[4,5,7,8,13,14]).
\end{split} \end{equation}

It is well known that the Bell polynomials (also called Tochard polynomials or exponential polynomials) are defined by
\begin{equation} \begin{split} \label{02}
e^{x(e^t-1)}=\sum_{n=0}^\infty B_{n}(x)\frac{t^n}{n!},\,\,\,\,(\textnormal{see}\,\,[1,5,9,10,12,16]).
\end{split} \end{equation}

From \eqref{01} and \eqref{02}, we note that
\begin{equation} \begin{split} \label{03}
B_{n}(x)&=e^{-x}\sum_{k=0}^\infty \frac{k^n}{k!}x^{k}\\
&=\sum_{k=0}^{n} x^{k} S_{2}(n,k),\,\,\,\,(n \geq 0),\,\,\,\,(\textnormal{see}\,\,[2,4,5]).
\end{split} \end{equation}
When $x=1$, $B_{n}=B_{n}(1)$ are called Bell numbers.\\

The (exponential) incomplete Bell polynomials (also called (exponential) partial Bell polynomials) are defined by the generating function 
\begin{equation} \begin{split} \label{04}
\frac{1}{k!}\Big(\sum_{m=1}^\infty x_{m}\frac{t^m}{m!}\Big)^{k}=\sum_{n=k}^\infty B_{n,k}(x_1,\cdots,x_{n-k+1})\frac{t^n}{n!},\,\,\,\,(k \geq 0),\,\,\,\,(\textnormal{see}\,\,[12,16]).
\end{split} \end{equation}

Thus, by \eqref{04}, we get
\begin{equation} \begin{split} \label{05}
B_{n,k}(x_1,\cdots,x_{n-k+1})&=\sum\frac{n!}{i_{1}!i_{2}!\cdots i_{n-k+1}!}\Big(\frac{x_{1}}{1!}\Big)^{i_{1}}\Big(\frac{x_{2}}{2!}\Big)^{i_{2}}\times \cdots\\
&\times \Big(\frac{x_{n-k+1}}{(n-k+1)!}\Big)^{i_{n-k+1}},
\end{split} \end{equation}
where the summation is over all integers $i_{1},\cdots,i_{n-k+1} \geq 0$ such that $i_{1}+i_{2}+\cdots+i_{n-k+1}=k$ and $i_{1}+2i_{2}+\cdots+(n-k+1)i_{n-k+1}=n$.\\

From \eqref{01} and \eqref{04}, we note that
\begin{equation} \begin{split} \label{06}
B_{n,k}\underbrace{(1,1,\cdots,1)}_{n-k+1-times}=S_{2}(n,k),\,\,\,\,(n,k \geq 0).
\end{split} \end{equation}

By \eqref{05}, we easily get
\begin{equation} \begin{split} \label{07}
B_{n,k}(\alpha x_{1},\alpha x_{2}, \cdots, \alpha x_{n-k+1})=\alpha^{k}B_{n,k}(x_{1},x_{2},\cdots,x_{n-k+1})
\end{split} \end{equation}
and
\begin{equation} \begin{split} \label{08}
B_{n,k}(\alpha x_{1},\alpha^{2} x_{2}, \cdots, \alpha^{n-k+1} x_{n-k+1})=\alpha^{n}B_{n,k}(x_{1},x_{2},\cdots,x_{n-k+1}),
\end{split} \end{equation}
where $\alpha \in \mathbb{R}$ $(\textnormal{see}\,\,[12,14])$.\\

From \eqref{04}, we easily note that
\begin{equation} \begin{split} \label{09}
&\sum_{n=k}^\infty B_{n,k}(x,1,0,0,\cdots,0)\frac{t^n}{n!}=\frac{1}{k!}\big(xt+\frac{t^2}{2}\big)^{k}\\
&=\frac{t^k}{k!}\sum_{n=0}^{k} {k \choose n}\Big(\frac{t}{2}\Big)^{n}x^{k-n}\\
&=\sum_{n=0}^{k}\frac{(n+k)!}{k!}\binom{k}{n}\frac{1}{2^n}x^{k-n}\frac{t^{n+k}}{(n+k)!},
\end{split} \end{equation}
and
\begin{equation} \begin{split} \label{10}
\sum_{n=k}^\infty B_{n,k}(x,1,0,0,\cdots,0)\frac{t^n}{n!}=\sum_{n=0}^\infty B_{n+k,k}(x,1,0,\cdots,0)\frac{t^{n+k}}{(n+k)!}.
\end{split} \end{equation}

By comparing the coefficients on both sides of \eqref{09} and \eqref{10}, we get
\begin{equation} \begin{split} \label{11}
B_{n+k,k}(x,1,0,\cdots,0)=\frac{(n+k)!}{k!}{k \choose n}\frac{1}{2^n}x^{k-n},\,\,\,(0 \leq n \leq k).
\end{split} \end{equation}

By replacing $n$ by $n-k$ in \eqref{11}, we get
\begin{equation} \begin{split} \label{12}
B_{n,k}(x,1,0,\cdots,0)=\frac{n!}{k!}{k \choose n-k}x^{2k-n}\Big(\frac{1}{2}\Big)^{n-k},\,\,\,(k \leq n \leq 2k).
\end{split} \end{equation}

The (exponential) complete Bell polynomials are defined by 
\begin{equation} \begin{split} \label{13}
\textnormal{exp}\Big(\sum_{i=1}^\infty x_{i}\frac{t^i}{i!}\Big)=\sum_{n=0}^\infty B_{n}(x_{1},x_{2},\cdots,x_{n})\frac{t^n}{n!}.
\end{split} \end{equation}

Then, by \eqref{04} and \eqref{13}, we get
\begin{equation} \begin{split} \label{14}
B_{n}(x_{1},x_{2},\cdots,x_{n})=\sum_{k=0}^{n}B_{n,k}(x_{1},x_{2},\cdots,x_{n-k+1})
\end{split} \end{equation}

From \eqref{03}, \eqref{06}, \eqref{07} and \eqref{14}, we have
\begin{equation} \begin{split} \label{15}
B_{n}(x,x,\cdots,x)&=\sum_{k=0}^{n}x^{k}B_{n,k}(1,1,\cdots,1)\\
&=\sum_{k=0}^{n} x^{k} S_{2}(n,k)=B_{n}(x),\,\,\,\,(n \geq 0).
\end{split} \end{equation}

It is known that the central factorial numbers of the second kind are given by 
\begin{equation} \begin{split} \label{16}
\frac{1}{k!}\big(e^{\frac{t}{2}}-e^{-\frac{t}{2}}\big)^{k}=\sum_{n=k}^\infty T(n,k)\frac{t^n}{n!},\,\,\,\,(\textnormal{see}\,\,[3,9,11]),
\end{split} \end{equation}
where $k \geq 0$.\\

From \eqref{16}, we can derive the following equation
\begin{equation} \begin{split} \label{17}
T(n,k)=\frac{1}{k!}\sum_{j=0}^{k} {k \choose j} (-1)^{k-j}\big(j-\frac{k}{2}\big)^{n},
\end{split} \end{equation}
where $n,k \in \mathbb{Z}$ with $n \geq k \geq 0$, $(\textnormal{see}\,\,[8,9,11])$.\\

In [11], the central Bell polynomials $B_{n}^{(c)}(x)$ are defined by
\begin{equation} \begin{split} \label{18}
B_{n}^{(c)}(x)=\sum_{k=0}^{n} T(n,k)x^{k},\,\,\,\,(n \geq 0).
\end{split} \end{equation}
When $x=1$, $B_{n}^{(c)}=B_{n}^{(c)}(1)$ are called the central Bell numbers.\\

From \eqref{18}, we can derive the generating function for the central Bell polynomials as follows:
\begin{equation} \begin{split} \label{19}
e^{x\big(e^{\frac{t}{2}}-e^{-\frac{t}{2}}\big)}=\sum_{n=0}^\infty B_{n}^{(c)}(x)\frac{t^n}{n!},\,\,\,\,(\textnormal{see}\,\,[9]).
\end{split} \end{equation}

Thus, by \eqref{19}, we have the following Dobinski-like
formula
\begin{equation} \begin{split} \label{20}
B_{n}^{(c)}(x)=\sum_{l=0}^\infty \sum_{j=0}^\infty {l+j \choose j}(-1)^{j}\frac{1}{(l+j)!}\Big(\frac{l}{2}-\frac{j}{2}\Big)^{n}x^{l+j},
\end{split} \end{equation}
where $n \geq 0$ $(\textnormal{see}\,\,[9])$.\\

Motivated by \eqref{04} and \eqref{13}, we introduce central complete and incomplete Bell polynomials and investigate some properties and identities for these polynomials.
Especially, we give explicit formulas for the central complete and incomplete Bell polynomials related to central factorial numbers of the second kind.

\section{\bf On central complete and incomplete Bell polynomials}
In view of \eqref{13}, we consider the \textit{central incomplete Bell polynomials} which are given by
\begin{equation} \begin{split} \label{21}
\frac{1}{k!}\Big(\sum_{m=1}^\infty \frac{1}{2^m}(x_{m}-(-1)^{m}x_{m})\frac{t^m}{m!}\Big)^{k}=\sum_{n=k}^\infty T_{n,k}(x_1,x_2,\cdots,x_{n-k+1})\frac{t^n}{n!},
\end{split} \end{equation}
where $k=0,1,2,3,\cdots$.\\

For $n,k \geq 0$ with $n-k\equiv 0$ (mod $2$), by \eqref{04} and \eqref{05}, we get
\begin{equation} \begin{split} \label{22}
T_{n,k}(x_1,x_2,\cdots,x_{n-k+1})&=\sum\frac{n!}{i_{1}!i_{2}!\cdots i_{n-k+1}!}\Big(\frac{x_{1}}{1!}\Big)^{i_{1}}\Big(\frac{0}{2\cdot 2!}\Big)^{i_{2}}\\
&\times(\frac{x_3}{2^2\cdot 3!}\Big)^{i_{3}}\cdots \Big(\frac{x_{n-k+1}}{2^{n-k}(n-k+1)!}\Big)^{i_{n-k+1}},
\end{split} \end{equation}
where the summation is over all integers $i_{1},i_{2},\cdots, i_{n-k+1} \geq 0$ such that $i_{1}+\cdots+i_{n-k+1}=k$ and $i_{1}+2i_{2}+\cdots+(n-k+1)i_{n-k+1}=n$.\\

From \eqref{05} and \eqref{22}, we note that
\begin{equation} \begin{split} \label{23}
T_{n,k}(x_{1},x_{2},\cdots,x_{n-k+1})=B_{n,k}\big(x_{1},0,\frac{x_{3}}{2^2},0,\cdots,\frac{x_{n-k+1}}{2^{n-k}}\big),
\end{split} \end{equation}
where $n,k \geq 0$ with $n-k \equiv0$ (mod $2$) and $n \geq k$.\\
Therefore, we obtain the following lemma.
\begin{lem} For $n,k \geq 0$ with $n \geq k$ and $n-k \equiv0$ (mod $2$), we have
\begin{equation*} \begin{split}
T_{n,k}(x_1,x_2,\cdots,x_{n-k+1})=B_{n,k}\big(x_{1},0,\frac{x_{3}}{2^2},0,\cdots,\frac{x_{n-k+1}}{2^{n-k}}\big).
\end{split} \end{equation*}
\end{lem}

\vspace{0.1in}

For $n,k \geq 0$ with $n \geq k$ and $n-k\equiv0$ (mod $2$), by \eqref{21}, we get
\begin{equation} \begin{split} \label{24}
\sum_{n=k}^\infty T_{n,k}(x,x^2,x^3,&\cdots,x^{n-k+1})\frac{t^n}{n!}=\frac{1}{k!}\Big(xt+\frac{x^3}{2^2}\frac{t^3}{3!}+\frac{x^5}{2^4}\frac{t^5}{5!}+\cdots\Big)^{k}\\
&=\frac{1}{k!}\Big(e^{\frac{x}{2}t}-e^{-\frac{x}{2}t}\Big)^{k}=\frac{1}{k!}e^{-\frac{kx}{2}t}\Big(e^{xt}-1\Big)^{k}\\
&=\frac{1}{k!}\sum_{l=0}^{k} {k \choose l}(-1)^{k-l}e^{(l-\frac{k}{2})xt}\\
&=\frac{1}{k!}\sum_{l=0}^{k} {k \choose l} (-1)^{k-l} \sum_{n=0}^\infty \big(l-\frac{k}{2}\big)^{n}x^{n}\frac{t^n}{n!}\\
&=\sum_{n=0}^\infty \Big(\frac{x^n}{k!}\sum_{l=0}^{k} {k \choose l} (-1)^{k-l} \big(l-\frac{k}{2}\big)^{n}\Big)\frac{t^n}{n!}.
\end{split} \end{equation}

Therefore, by comparing the coefficients on both sides of \eqref{24}, we obtain the following theorem.
\begin{thm} For $n,k \geq 0$ with $n-k \equiv 0$ (mod $2$), we have
\begin{equation} \begin{split} \label{25}
 \frac{x^n}{k!}\sum_{l=0}^{k} {k \choose l} (-1)^{k-l} \big(l-\frac{k}{2}\big)^{n}=
 \left\{ \begin{array}{lcr}
T_{n,k}(x,x^2,\cdots,x^{n-k+1}),& \textnormal{if}\, n \geq k, \\
 0,  \quad   \quad & \textnormal{if}\, n<k.  \end{array} \right.
\end{split} \end{equation}
In particular
\begin{equation} \begin{split} \label{26}
 \frac{1}{k!}\sum_{l=0}^{k} {k \choose l} (-1)^{k-l} \big(l-\frac{k}{2}\big)^{n}=
 \left\{ \begin{array}{lcr}
T_{n,k}(1,1,\cdots,1),& \textnormal{if}\, n \geq k, \\
 0,  \quad   \quad & \textnormal{if}\, n<k.  \end{array} \right.
\end{split} \end{equation}
\end{thm}

For $n,k \geq 0$ with $n-k \equiv 0$ (mod $2$) and $n \geq k$, by \eqref{17} and \eqref{26}, we get
\begin{equation} \begin{split} \label{27}
T_{n,k}(1,1,\cdots,1)=T(n,k).
\end{split} \end{equation}

Therefore, by \eqref{25} and \eqref{26}, we obtain the following corollary
\begin{cor} For $n,k \geq 0$ with $n-k \equiv 0$ (mod $2$), $n \geq k$, we have
\begin{equation*} \begin{split}
&T_{n,k}(x,x^2,\cdots,x^{n-k+1})=x^{n}T_{n,k}(1,1,\cdots,1)
\end{split} \end{equation*}
and
\begin{equation*} \begin{split}
&T_{n,k}(1,1,\cdots,1)=T(n,k)=B_{n,k}\big(1,0,\frac{1}{2^2},\cdots,\frac{1}{2^{n-k}}\big)\\
&=\sum \frac{n!}{i_{1}!i_{3}!\cdots i_{n-k+1}!}\Big(\frac{1}{1!}\Big)^{i_{1}}\Big(\frac{1}{2^{2}3!}\Big)^{i_{3}}\cdots\Big(\frac{1}{2^{n-k}(n-k+1)!}\Big)^{i_{n-k+1}},
\end{split} \end{equation*}
where $i_{1}+i_{3}+\cdots+i_{n-k+1}=k$ and $i_{1}+3i_{3}+\cdots +(n-k+1)i_{n-k+1}=n$.
\end{cor}

For $n,k \geq 0$ with $n\geq k$ and $n-k \equiv0$ (mod $2$), we observe that
\begin{equation} \begin{split} \label{28}
\sum_{n=k}^\infty T_{n,k}(x,1,0,0,\cdots,0)\frac{t^n}{n!}=\frac{1}{k!}(xt)^{k}.
\end{split} \end{equation}

Thus we have
\begin{equation*} \begin{split}
T_{n,k}(x,1,0,\cdots,0)=x^{k}{0 \choose n-k}.
\end{split} \end{equation*}

For $n,k \geq 0$ with $n-k \equiv 0$ (mod $2$), $n\geq k$, by \eqref{22}, we get
\begin{equation} \begin{split} \label{29}
T_{n,k}(x_1,x_2,\cdots,x_{n-k+1})&=\sum\frac{n!}{i_{1}!i_{3}!\cdots i_{n-k+1}!}\Big(\frac{x_{1}}{1!}\Big)^{i_{1}}(\frac{x_3}{2^2\cdot 3!}\Big)^{i_{3}}\\
&\times \cdots \times \Big(\frac{x_{n-k+1}}{2^{n-k}(n-k+1)!}\Big)^{i_{n-k+1}},
\end{split} \end{equation}
where the summation is over all integers $i_1,i_2,\cdots,i_{n-k+1}\geq 0$ such that $i_{1}+i_{3}+\cdots+i_{n-k+1}=k$ and $i_{1}+3i_{3}+\cdots+(n-k+1)i_{n-k+1}=n$.\\

By \eqref{29}, we easily get
\begin{equation} \begin{split} \label{30}
T_{n,k}(x,x,\cdots,x)=x^{k}T_{n,k}(1,1,\cdots,1)
\end{split} \end{equation}
and
\begin{equation*} \begin{split}
T_{n,k}(\alpha x_{1},\alpha x_{2},\cdots,\alpha x_{n-k+1})=\alpha^{k} T_{n,k}(x_{1},x_{2},\cdots,x_{n-k+1}),
\end{split} \end{equation*}
where $n,k \geq 0$ with $n-k \equiv 0$ (mod $2$) and $n \geq k$.\\

Now, we observe that
\begin{equation} \begin{split} \label{31}
\textnormal{exp}&\Big(x\sum_{i=1}^\infty \big(\frac{1}{2}\big)^{i}(x_{i}-(-1)^{i}x_{i})\frac{t^i}{i!}\Big)\\
&=\sum_{k=0}^\infty x^{k} \frac{1}{k!}\Big(\sum_{i=1}^\infty \big(\frac{1}{2}\big)^{i}(x_{i}-(-1)^{i}x_{i})\frac{t^i}{i!}\Big)^{k}\\
&=1+\sum_{k=1}^\infty x^{k} \frac{1}{k!}\Big(\sum_{i=1}^\infty \big(\frac{1}{2}\big)^{i}(x_{i}-(-1)^{i}x_{i})\frac{t^i}{i!}\Big)^{k}\\
&=1+\sum_{k=1}^\infty x^{k}\sum_{n=k}^\infty T_{n,k}(x_{1},x_{2},\cdots,x_{n-k+1})\frac{t^n}{n!}\\
&=1+\sum_{n=1}^\infty \Big(\sum_{k=1}^{n} x^{k} T_{n,k}(x_{1},x_{2},\cdots,x_{n-k+1})\Big)\frac{t^n}{n!}.
\end{split} \end{equation}

In view of \eqref{13}, we define the \textit{central complete Bell polynomials} by
\begin{equation} \begin{split} \label{32}
\textnormal{exp}\Big(x\sum_{i=1}^\infty \big(\frac{1}{2}\big)^{i}(x_{i}-(-1)^{i}x_{i})\frac{t^i}{i!}\Big)=\sum_{n=0}^\infty B_{n}^{(c)}(x|x_{1},x_{2},\cdots,x_{n})\frac{t^n}{n!}.
\end{split} \end{equation}

Thus, by \eqref{31} and \eqref{32}, we get
\begin{equation} \begin{split} \label{33}
B_{n}^{(c)}(x|x_{1},x_{2},\cdots,x_{n})=\sum_{k=0}^{n} x^{k}T_{n,k}(x_{1},x_{2},\cdots,x_{n-k+1}).
\end{split} \end{equation}
When $x=1$, $B_{n}^{(c)}(1|x_{1},x_{2},\cdots,x_{n})=B_{n}^{(c)}(x_{1},x_{2},\cdots,x_{n})$ are called the \textit{central complete Bell numbers}.\\

For $n \geq 0$, we have
\begin{equation} \begin{split} \label{34}
B_{n}^{(c)}(x_{1},x_{2},\cdots,x_{n})=\sum_{k=0}^{n} T_{n,k}(x_{1},x_{2},\cdots,x_{n-k+1})
\end{split} \end{equation}
and
\begin{equation*} \begin{split}
B_{0}^{(c)}(x_{1},x_{2},\cdots,x_{n})=1.
\end{split} \end{equation*}

By \eqref{18} and \eqref{33}, we get
\begin{equation} \begin{split} \label{35}
B_{n}^{(c)}(1,1,\cdots,1)=\sum_{k=0}^{n} T_{n,k}(1,1,\cdots,1)=\sum_{k=0}^{n} T(n,k)=B_{n}^{(c)},
\end{split} \end{equation}
and
\begin{equation} \begin{split} \label{36}
B_{n}^{(c)}(x|1,1,\cdots,1)=\sum_{k=0}^{n} x^{k}T_{n,k}(1,1,\cdots,1)=\sum_{k=0}^{n} x^{k}T(n,k)=B_{n}^{(c)}(x).
\end{split} \end{equation}

From \eqref{31}, we note that
\begin{equation*} \begin{split}
\textnormal{exp}\Big(\sum_{i=1}^\infty &\big(\frac{1}{2}\big)^{i}(x_{i}-(-1)^{i}x_{i})\frac{t^i}{i!}\Big)\\
&=1+\sum_{n=1}^\infty \frac{1}{n!}\Big(\sum_{i=1}^\infty \big(\frac{1}{2}\big)^{i}(x_{i}-(-1)^{i}x_{i})\frac{t^i}{i!}\Big)^{n}\\
&=1+\frac{1}{1!}\sum_{i=1}^\infty \big(\frac{1}{2}\big)^{i}(x_{i}-(-1)^{i}x_{i})\frac{t^i}{i!}+\frac{1}{2!}\Big(\sum_{i=1}^\infty \big(\frac{1}{2}\big)^{i}(x_{i}-(-1)^{i}\\
&\times x_{i})\frac{t^i}{i!}\Big)^{2}+\frac{1}{3!}\Big(\sum_{i=1}^\infty \big(\frac{1}{2}\big)^{i}(x_{i}-(-1)^{i}x_{i})\frac{t^i}{i!}\Big)^{3}+\cdots\\
\end{split} \end{equation*}
\begin{equation} \begin{split} \label{37}
&=1+\frac{1}{1!}x_{1}t+\frac{1}{2!}x_{1}^{2}t^2+\Big(\frac{1}{3!2^2}x_{3}+\frac{x_{1}^3}{3!}\Big)t^3+\cdots\\
&=\sum_{n=0}^\infty \Big(\sum_{m_{1}+2m_{2}+\cdots+nm_{n}=n} \frac{n!}{m_{1}!m_{2}!\cdots m_{n}!}\Big(\frac{x_{1}}{1!}\Big)^{m_{1}}\Big(\frac{0}{2!}\Big)^{m_{2}}\\
&\times \Big(\frac{x_{3}}{3!2^2}\Big)^{m_{3}}\cdots \Big(\frac{x_{n}\big(1-(-1)^{n}\big)}{n!2^n}\Big)^{m_{n}}\Big)\frac{t^n}{n!}.
\end{split} \end{equation}

\vspace{0.1in}

Now, for $n \in \mathbb{N}$ with $n \equiv 1$ (mod $2$), by \eqref{32}, \eqref{34} and \eqref{37}, we get
\begin{equation} \begin{split} \label{38}
&B_{n}^{(c)}(x_{1},x_{2},\cdots,x_{n})=\sum_{k=0}^{n} T_{n,k}(x_{1},x_{2},\cdots,x_{n-k+1})\\\
&=\sum_{m_{1}+3m_{3}+\cdots+nm_{n}=n}\frac{n!}{m_{1}!m_{3}!\cdots m_{n}!}\Big(\frac{x_{1}}{1!}\Big)^{m_{1}}
\Big(\frac{x_{3}}{3!2^2}\Big)^{m_{3}} \cdots\Big(\frac{x_{n}}{n!2^{n-1}}\Big)^{m_{n}}.
\end{split} \end{equation}

Therefore, by \eqref{38}, we obtain the following theorem
\begin{thm} For $n \in \mathbb{N}$ with $n \equiv 1$ (mod $2$), we have
\begin{equation*} \begin{split}
B_{n}^{(c)}(x_{1},x_{2},\cdots,x_{n})=&=\sum_{m_{1}+3m_{3}+\cdots+nm_{n}=n}\frac{n!}{m_{1}!m_{3}!\cdots m_{n}!}\Big(\frac{x_{1}}{1!}\Big)^{m_{1}}\\
&\times \Big(\frac{x_{3}}{3!2^2}\Big)^{m_{3}} \cdots\Big(\frac{x_{n}}{n!2^{n-1}}\Big)^{m_{n}}.
\end{split} \end{equation*}
\end{thm}

We note that
\begin{equation} \begin{split} \label{39}
\textnormal{exp}\Big(x\sum_{i=1}^\infty \big(\frac{1}{2}\big)^{i}\big(1-(-1)^{i}\big)\frac{t^i}{i!}\Big)
&=1+\sum_{k=1}^\infty \frac{x^k}{k!}\Big(\sum_{n=k}^\infty \big(\frac{1}{2}\big)^{i}\big(1-(-1)^{i}\big)\frac{t^i}{i!}\Big)^{k}\\
&=1+\sum_{k=1}^{\infty}x^k\sum_{n=k}^{\infty}T_{n,k}(1,1,\cdots,1)\frac{t^n}{n!}\\
&=1+\sum_{n=1}^\infty \Big(\sum_{k=1}^{n} x^{k}T_{n,k}(1,1,\cdots,1)\Big)\frac{t^n}{n!}.
\end{split} \end{equation}

On the other hand, from \eqref{19} we have
\begin{equation} \begin{split} \label{40}
\textnormal{exp}\Big(x\sum_{i=1}^\infty \big(\frac{1}{2}\big)^{i}\big(1-(-1)^{i}\big)\frac{t^i}{i!}\Big)&=\textnormal{exp}\Big(x\big(t+\frac{1}{2^2}t^{3}+\frac{1}{2^4}t^5+\cdots\big)\Big)\\
&=\textnormal{exp}\Big(x\big(e^{\frac{t}{2}}-e^{-\frac{t}{2}}\big)\Big)=\sum_{n=0}^\infty B_{n}^{(c)}(x)\frac{t^n}{n!}.
\end{split} \end{equation}

Therefore, by \eqref{39} and \eqref{40}, we obtain the following theorem.
\begin{thm} For $n,k\geq0$ with $n \geq k$, we have
\begin{equation*} \begin{split}
\sum_{k=0}^{n} x^{k}T_{n,k}(1,1,\cdots,1)=B_{n}^{(c)}(x).
\end{split} \end{equation*}
\end{thm}

From Theorem 2.5, we note that
\begin{equation} \begin{split} \label{41}
\sum_{k=0}^{n} x^{k}T_{n,k}(1,1,\cdots,1)=\sum_{k=0}^{n} T_{n,k}(x,x,\cdots,x)=B_{n}^{(c)}(x,x,\cdots,x).
\end{split} \end{equation}

Therefore, by Theorem 2.5 and \eqref{41}, we obtain the following corollary.
\begin{cor} For $n \geq 0$, we have
\begin{equation*} \begin{split}
B_{n}^{(c)}(x,x,\cdots,x)=B_{n}^{(c)}(x).
\end{split} \end{equation*}
\end{cor}

It is known that the Stirling numbers of the first kind are given by the generating function 
\begin{equation} \begin{split} \label{42}
\frac{1}{k!}\big(\log (1+t)\big)^{k}=\sum_{n=k}^\infty S_{1}(n,k)\frac{t^n}{n!},\,\,\,\,(k \geq 0),\,\,\,\,(\textnormal{see}\,\,[4,6]).
\end{split} \end{equation}

By \eqref{42}, we easily get
\begin{equation} \begin{split} \label{43}
\frac{1}{k!}\log\Big(1+\frac{x}{1-\frac{x}{2}}\Big)&=\sum_{l=k}^\infty S_{1}(l,k)\frac{1}{l!}\Big(\frac{x}{1-\frac{x}{2}}\Big)^{l}\\
&=\sum_{l=k}^\infty S_{1}(l,k)\frac{x^l}{l!}\big(1-\frac{x}{2}\big)^{-l}\\
&=\sum_{l=k}^\infty \frac{1}{l!}S_{1}(l,k)\sum_{n=l}^\infty {n-1 \choose l-1} \big(\frac{1}{2}\big)^{n-l}x^{n}\\
&=\sum_{n=k}^\infty \Big(\sum_{l=k}^{n} \frac{1}{l!}S_{1}(l,k){n-1 \choose l-1} \big(\frac{1}{2}\big)^{n-l}\Big)x^n.
\end{split} \end{equation}

From \eqref{21} and \eqref{43}, we can derive the following equation.
\begin{equation} \begin{split} \label{44}
\sum_{n=k}^\infty& T_{n,k}(0!,1!,2!,\cdots,(n-k)!\Big)\frac{t^n}{n!}\\
&=\frac{1}{k!}\Big(t+\big(\frac{1}{2}\big)^{2}\frac{t^3}{3}+\big(\frac{1}{2}\big)^{4}\frac{t^5}{5}+\big(\frac{1}{2}\big)^{6}\frac{t^7}{7}+\cdots\Big)^{k}\\
&=\frac{1}{k!}\Big(\log \big(1+\frac{t}{2}\big)-\log \big(1-\frac{t}{2}\big)\Big)^{k}=\frac{1}{k!}\bigg(\log \Big(\frac{1+\frac{t}{2}}{1-\frac{t}{2}}\Big)\bigg)^{k}\\
&=\frac{1}{k!}\Big(\log\big(1+\frac{t}{1-\frac{t}{2}}\big)\Big)^{k}=\sum_{n=k}^\infty \Big(\sum_{l=k}^{n} \frac{S_{1}(l,k)}{l!}{n-1 \choose l-1}\big(\frac{1}{2}\big)^{n-l}\Big)t^n.
\end{split} \end{equation}

By comparing the coefficients on both sides of \eqref{44}, we obtain the following theorem.
\begin{thm} For $n,k \geq 0$ with $n \geq k$, we have
\begin{equation*} \begin{split}
T_{n,k}\big(0!,1!,2!,\cdots,(n-k)!\big)=n!\sum_{l=k}^{n} \frac{S_{1}(l,k)}{l!}{n-1 \choose l-1} \big(\frac{1}{2}\big)^{n-l}.
\end{split} \end{equation*}
\end{thm}

\end{document}